\documentclass[12pt, 14paper,reqno]{amsart}
\usepackage{amsmath,amsfonts,amssymb, fancyhdr, lipsum}

\theoremstyle{plain}
\newtheorem{theorem}{Theorem}
\newtheorem{theorem*}{Theorem}
\newtheorem{lemma}{Lemma}
\newtheorem{remark}{Remark}
\newtheorem{definition}{Definition}

\newtheorem{proposition*}{Proposition}

\theoremstyle{remark}
\numberwithin{equation}{section}
\numberwithin{lemma}{section}
\numberwithin{remark}{section}
\numberwithin{theorem}{section}
\numberwithin{definition}{section}

\usepackage{amsmath}
\usepackage{amsfonts}   
\usepackage{amssymb}
\usepackage{amssymb, amsmath, amsthm}
\usepackage[breaklinks]{hyperref}
\newtheorem{example}{Example}

\usepackage{graphicx}

\begin{document}
\title[General solution of the Diophantine equation $M_p^{x} +  (M_q+1)^{y}=(lz)^2$] {General solution of the Diophantine equation $M_p^{x} +  (M_q+1)^{y}=(lz)^2$}

\author {Ghosh, Arkabrata}
\address{Arkabrata Ghosh, Department of Mathematics, Central Michigan University, Mount Pleasant, Michigan 48858, USA}
\email{arka2686@gmail.com}
\keywords{Diophantine equation, Mersenne Primes, Integer solution}
\subjclass[2020] {11D61, 11D72, 11A41}

\maketitle

\begin{abstract}
    In this article, I study and solve the exponential Diophantine equation $M_p^{x} + (M_q + 1)^{y}= (lz)^2$ where $M_p$ and $M_q$ are Mersenne primes, $l$ is a prime number, and $x$,$y$  and $z$ are non-negative integers. Several illustrations are presented as well as cases where no solution of the given Diophantine equation is present.
\end{abstract}

\section{Introduction}
 
    The Diophantine equation is one of the most attractive and exciting categories of problems in number theory. over the years, a number of researchers have been studying the Diophantine equation of the form $a^x + b^y=z^2$. This includes Aggarwal, Burshtein, Sroysang, Rabago, among others(\cite{A20}, \cite{AN20}, \cite{AS17}). Some of them have studied these equations in relation to Mersenne primes( see definition \ref{def:mer}). The primary focus of their work is on the case where one of the bases $a$ and $b$ is a Mersenne prime. Sroysang ~\cite{S13} proved that the solutions of $3x + 2y = z^2$ are $(0, 1, 2)$; $(3,0,3)$ and $(2, 4,5)$. Asthana and Singh ~\cite{AS17} proved that $3x + 13y = z^2$ has exactly four non-negative integer solutions, and these are $(1, 0, 2)$, $(1, 1, 4)$; $(3, 2, 14)$ and $(5, 1, 6)$. Rabago ~\cite{R13} proved that the triples $(4, 1, 10)$ and $(1, 0, 2)$ are the only solutions to the Diophantine equation $3^x + 19^y = z^2$, and that $(2,1, 10)$ and $(1, 0, 2)$ are the only two solutions to $3^x + 91^y = z^2$. Sroysang ~\cite{S13} also showed that the $7^x + 8^y = z^2$ has the only solution $(x, y, z) = (0, 1, 3)$.  Chotchaisthit~\cite{C13} aimed to study $px + (p + 1)y = z^2$ in the set of non-negative integers and where $p$ is a Mersenne prime. 
    \newpage

    In this article, I have found a general solution of the exponential Diophantine equation $M_p^{x} +  (M_q+1)^{y}=(lz)^2$ where $M_p$ and $M_q$ are Mersenne Primes, $x$, $y$ and $z$ are non-negative integers and $l$ is a prime number. Methods of modular arithmetic and factorization of polynomials are used in proving the results of this article.

\section{Main Results}

The following definition and lemmas are needed for this article.
\begin{definition}
\label{def:mer}
      A Mersenne prime is a prime number of the form $2^p-1$ where $p$ is a prime number and is denoted by $M_p$. 
\end{definition}

  \begin{lemma}
  \label{le:Mersenne}
      All Mersenne Primes are congruent to $3 \pmod 4$
  \end{lemma}

  \begin{proof}
      As any Mersenne Prime is of the form $2^p-1$, we can clearly say that $ p \geq 2$. Now as $ p\geq 2$, then $2^p \equiv 0 \pmod 4$ and hence, $ 2^p-1 \equiv 3 \pmod 4$.
  \end{proof}

\begin{lemma}
{(Mihailescu's theorem)}(see ~\cite{M04}) The only solution to the the Diophantine equation $a^x-b^y=1$ is $a=3$, $b=2$, $x=2$ and $y=2$ where min$\{a,b,x,y\} > 1$.
 \label{le:mihe}   
\end{lemma}
At first, we consider the case when $l=2$. Hence, the following is the first main theorem of this paper.
  \begin{theorem}
  \label{thm1}
     Every non-negative integer solution to the equation
      $M_p^x + (M_q+ 1)^y= (2z)^2$ is the following tuple $(M_p, M_q,x,y,z) =(3, M_q,1,0,1)$.
  \end{theorem}

  \begin{proof}
      Let us first consider the case when one of the exponent $x$ and $y $ is zero. \\
      Case-I At first, if we assume $x=0$, then we get the following equation
      \begin{equation}
          2^{qy}= 4z^2-1.
          \label{eq:2.1}
      \end{equation}
      Subcase-(a) If $y=0$, then from the equation \ref{eq:2.1}, we get $4z^2=2$ which is a contradiction.\\
      Subcase-(b) If $y=1$, then from the equation \ref{eq:2.1}, we get that $(2z)^2- 2^{q}=1$. From the lemma \ref{le:mihe}, solution of this equation is only possible if $z= \frac{3}{2}$ which is a contradiction. \\
      Subcase-(c) If $y > 1$, then the equation \ref{eq:2.1} can be written as $(2z)^2- 2^{qy}=1$. By the lemma \ref{le:mihe}, we must have $qy=3$. Now as $q$ is a prime number, we get $q=3$ and $y=1$ which is a contradiction to our assumption.\\
     Case-II Now we assume $y=0$. Then we get the following equation
      \begin{equation}
      \label{eq:2.2}
          M_p^{x} + 1 = (2z)^2.
      \end{equation}
   Subcase-(a) If $x=0$, then from the equation \ref{eq:2.2}, we get $4z^2=2$ which is a contradiction.\\
   Subcase-(b) If $x=1$, then the equation \ref{eq:2.2} can be written as $(2z)^2= 2^{p}$. Now let $Z=2z$ and $Z=2^{a}$. Then we get $2^{2a}=2^{p}$ which in turn gives $p=2a$. As $p$ is a prime number, we get $a=1$ and $p=2$. Hence $Z=2$ and finally, $z=1$. Hence $(M_p, M_q, x,y,z)=(3, M_q, 1.0,1)$.\\
   Subcase-(c) if $x >1$, then from the equation \ref{eq:2.2}, we get that $(2z)^2-(2^p-1)^x=1$. By the lemma \ref{le:mihe}, we get that $2^p=3$ which is a contradiction.\\
Case-III. Now we consider the case when $\{x,y\} \geq 1$. From the lemma \ref{le:Mersenne}, we know that $M_p \equiv 3 \pmod 4$ and $(M_q + 1) \equiv 0 \pmod 4$. Hence, 
$$
(M_p^{x} + (M_q+1)^{y}) \equiv \left\{\begin{array}{cc}
3 \pmod 4, &  \mbox{ $x$ is odd} \\
 1 \pmod 4, &  \mbox{ $x$ is even}
\end{array}\right.
$$

Now as $4z^2 \equiv 0 \pmod 4$, the equation has no solution when $\{x,y\} \geq 1$.
  \end{proof}

  Now we will consider the case when $l$ is an odd prime. Then we get the following theorem
  \begin{theorem}
  \label{thm2}
      Every non-negative integer solution to the equation $M_p^x + (M_q+ 1)^y= (lz)^2$ takes one of the following form:\\
      (a)$(M_p, M_q,x,y,z)$= $(M_p, 7, 0,1,1)$ \\
      (b) $(M_p, M_q,x,y,z)$ = $\bigg(M_p, M_q, 2, \frac{p+2}{q}, \frac{2^p+1}{l}\bigg)$.
      \end{theorem}

\begin{proof}
   Let us first assume that one of the exponents $x$ and $y$ is zero.\\
   Case-I At first, we assume that $x=0$.Then we get the equation
   \begin{equation}
       \label{eq:2.2.1}
       1 + (M_q+1)^y=(lz)^2.
   \end{equation}

\noindent Subcase-(a) If $y=0$, then from the equation \ref{eq:2.2.1}, we get $2=(lz)^2$ which is a contradiction.\\
Subcase-(b) If $y=1$, then from the equation \ref{eq:2.2.1}, we get that $(lz)^2- 2^q=1$. Hence from the lemma \ref{le:mihe}, we can conclude that $l=3$, $z=1$ and $q=3$. Hence, $(M_p, M_q, x,y,z)=(M_p, 7, 0,1,1) $ is the only solution when $l=3$.\\
Subcase-(c) If $ y > 1$, then from the equation \ref{eq:2.2.1}, we get $(lz)^2- 2^{qy}=1$. Again, from the lemma \ref{le:mihe}, we get that $qy=3$. By using the primality of $q$, we get that $q=3$ and $y=1$ which is a contradiction.

Case-II Now we assume $y=0$. Then we get the equation
   \begin{equation}
       \label{eq:2.2.2}
     (M_p)^x + 1=(lz)^2.
   \end{equation}
 Subcase-(a)  If $x=0$, then from the equation \ref{eq:2.2.2}, we get $2=(lz)^2$ which is a contradiction.\\
Subcase-(b) If $x=1$, then from the equation \ref{eq:2.2.2}, we get
$2^p=(lz)^2$. Now as $2^p\not \equiv 0 \pmod l$ and $(lz)^2 \equiv 0 \pmod l$, so the above equation has no solution.\\
Subcase-(c) If $x > 1$, then from the equation \ref{eq:2.2.2}, we get that $(lz)^2-(M_p)^x=1$. By the lemma \ref{le:mihe}, we get that $2^p=3$ which is a contradiction.

Case-III Now we consider the case when $ \{x,y\} \geq 1$. Now as we know $(lz)^2 \equiv 1 \pmod 4$ when $z$ is odd and 

$$
(M_p^{x} + (M_q+1)^{y}) \equiv \left\{\begin{array}{cc}
3 \pmod 4, &  \mbox{ $x$ is odd} \\
 1 \pmod 4, &  \mbox{ $x$ is even}
 \end{array}\right.
$$

. From the above, we can conclude that the equation $M_p^x + (M_q+ 1)^y= (lz)^2$ has a solution only if $x$ is even and $z$ is odd. Thus there exists a positive integer $k$ such that $x=2k$ and we get the equation $M_p^{2k} + 2^{qy}=(lz)^2$. This equation can be written as
\begin{equation}
    \label{eq: 2.2.3}
    (lz + M_p^k) (lz-M_p^k)=2^{qy}.
\end{equation}

There exist two non-negative integers $\alpha$ and $\beta$ with $\alpha > \beta $ such that 
$\alpha + \beta = qy$. Then the equation \ref{eq: 2.2.3} can be written as 
\begin{equation}
    \label{eq:2.2:4}
    (lz + M_p^k) (lz-M_p^k)=2^{\alpha + \beta}.
\end{equation}
we claim that gcd$( lz + M_p^k, lz-M_p^k) \neq 1$. Suppose our assumption is wrong. Then $( lz + M_p^k, lz-M_p^k)=1$. Now from the equation \ref{eq:2.2:4}, we can say that $lz-M_p^k=1$. We know from the lemma \ref{le:Mersenne} that $M_p \equiv 3 \pmod 4$ and hence $M_p^k \equiv 1,3 \pmod 4$.
As $l$ is an odd prime and $z$ is also odd, then $lz \equiv 1, 3 \pmod 4$. So, $lz-M_p^k \equiv 0, 2 \pmod 4$ which is contradiction to the fact $lz-M_p^k \equiv 1 \pmod 4$. Now as gcd$( lz + M_p^k, lz-M_p^k) \neq 1$, we take $lz+ M_p^k= 2^{\alpha}$ and  $ lz- M_p^k= 2^{\beta}$. it implies that $2M_p^k= 2^{\beta}(2^{\alpha - \beta} -1)$  and by comparing odd and even parts, we get the system of equations:

$$
\left\{\begin{array}{cc}
2^{\beta}=2.  \\
2^{\alpha-\beta}-1= M_p^{k}. 
 \end{array}\right.
$$

From the above equation, we know that $\beta =1$ and hence,

\begin{equation}
\label{eq:2.2.5}
    1= 2^{\alpha-1}-M_p^{k}.
\end{equation}

By using the lemma \ref{le:mihe}, we can say that the equation \ref{eq:2.2.5} has no solution if $\alpha > 2$  and $  k > 1$. If $\alpha=2 $, then from the equation $\ref{eq:2.2.5}$, we get that $M_p^{k}=1$. This gives $k=0$ which is a contradiction to that $k$ being a positive integer. So the only possibility is $k=1$ and hence $x=2$. Now putting the value $k=1$ in the equation \ref{eq:2.2.5}, we get $M_p + 1 = 2^{\alpha -1}$ or equivalently, $ \alpha= p+1$. Putting the values $\alpha= p+1$ and $\beta =1$ in the relation $\alpha \beta = qy$, we get $y= \frac{p+2}{q}$. Now putting the values $x$ and $y$ in the equation $M_p^x + (M_q+ 1)^y= (lz)^2$, we get that $z= \frac{2^p +1}{l}$.

\end{proof}
\begin{remark}
    From the theorem \ref{thm2}, we can say that the equation $M_p^x + (M_q+ 1)^y= (lz)^2$ has a positive integer solution $\bigg(M_p, M_q, 2, \frac{p+2}{q}, \frac{2^p+1}{l}\bigg)$. It means that given a Mersenne prime $M_p$, the solution can only be found if $q$ divides $(p+2)$ and also $l $ divides $(2^p + 1)$ where $l$ is an odd prime. We need to be careful about $q$ and pick only those $q$ for which $M_q$ is a Mersenne prime.
\end{remark}

      \begin{example}
          Find all possible positive integer solutions of the equation $8191^x + (M_q+1)^y= (3z)^2$ where $M_q$ is a Mersenne prime and $q$ is a prime number.
      \end{example}
          Solution:
              Here $M_p=8191 $ implies $p=13$. From the theorem \ref{thm2}, we get the condition that $q$ divides $p+2 =15$ which in turn gives $q=3$ or $q=5$. If $q=3$, we have $M_q= 7$, $y=5 $ and when $q=5$, we have $M_5=2^5-1=31$ and $y=3$. Observe that both $M_3$ and $M_5$ are Mersenne primes. Also for both the cases, $l=3$ and hence $z= \frac{2^{13}+1}{3}=2731$. Hence the solution set are  $(M_p, M_q,x,y,z)= ( 8191,7,2,5,2731) $ and $(M_p, M_q,x,y,z)= ( 8191,31,2,3,2731) $.

    \begin{example}
         Find the  positive integer solutions of the equation $7^{x}  + 4^y = (7z)^{2}$
    \end{example}

    Solution: Here $M_p=7$ where $p=3$ and $M_q=3$ where $q=2$. Hence, from the theorem \ref{thm2}, the solution set is $(x,y,z)=( 2, \frac{p+2}{q}, \frac{2^p+1}{l})$ if $q$ divides $(p+2)$ and $l$ divides $2^{p} + 1$. Now as $l=7$  does not divide $2^{p}+ 1=9$, the equation has no solution in positive integers.

    \begin{example}
        Find the positive integer solutions of the equation $3^{x}  + 8^y = (5z)^{2}$.
        \end{example}          
          Solution: Here $M_p=3$ where $p=2$ and $M_q=7$ where $q=3$. Hence, from the theorem \ref{thm2}, the solution set is $(x,y,z)=( 2, \frac{p+2}{q}, \frac{2^p+1}{l})$ if $q$ divides $(p+2)$ and $l$ divides $2^{p} + 1$. As $q=3$  does not divide $p+2= 4$, the equation has no solution in positive integers.
   
\section{Conclusion and future work}
  In this article, using the modular arithmetic method, with the help of Mihailescu's theorem \ref{le:mihe} and using the fact that every Mersenne prime is of the form $4k+3$, we have been able to show the complete list of positive integer solutions of the Diophantine equation $M_p^{x}  + (M_q+1)^y=(lz)^2$ where $l$ is a prime. 

  The following table represents the solution of the Diophantine equation  $M_p^{x}  + (M_q+1)^y=(lz)^2$ for the first couple of Mersenee Primes:

\begin{table}[ht]
    \centering

\begin{tabular}{|c| c| c| c| c |c |c |c|} 
 \hline
 $M_p$ & $p$ & $p+2$ & $q$ & $M_q$ & $2^{p}+1$ & $l$ & $(x,y,z)$ \\[0.5ex] 
 \hline\hline
 3 & 2 & 4 & 2 & 3 & 5 & 5 & (2,2,1) \\ 
 \hline
 7 & 3 & 5 & 5 & 7 & 9 & 3 & (2,1,3)\\
 \hline
 31 & 5 & 7 & 7 & 31 & 33 & 3 & (2,1,11) \\
 \hline
 31 & 5 & 7 & 7 & 31 & 33 & 11 & (2,1,3)\\
 \hline
 127 & 7 & 9 & 3 & 7 & 129 & 3 & (2,1,43)\\
 \hline
 127 & 7 & 9 & 3 & 7 & 129 & 43 & (2,1,3)\\ [1ex] 
 \hline
\end{tabular}
 \caption{Some possible integer solution of the equation $M_p^{x}  + (M_q+1)^y=(lz)^2$}
    \label{tab:my_label}
\end{table}

The next table represents some particular cases of the Diophantine equation $M_p^{x}+ (M_q+1)^y=(lz)^2$ where no solutions can be obtained. The un-solvability of these equations is due to two main reasons namely $q$ does not divide $(p+2)$ or $l$ does not divide $(2^p+1)$.
\begin{table}[ht]
    \centering

\begin{tabular}{|c| c| c| c| c |c |c |c|} 
 \hline
 $M_p$ & $p$ & $p+2$ & $q$ & $M_q$ & $2^{p}+1$ & $l$ & $M_p^{x}+ (M_q+1)^y=(lz)^2$ \\[0.5ex] 
 \hline\hline
 3 & 2 & 4 & 5 & 31 & 5 & 3 & $3^x + 32^y=(3z)^2$ \\ 
 \hline
 7 & 3 & 5 & 7 & 127 & 9 & 5 & $7^x + 128^y=(5z)^2$ \\
 \hline
 31 & 5 & 7 & 3 & 7 & 33 & 7 & $31^x + 8^y=(7z)^2$  \\
 \hline
 127 & 7 & 9 & 5 & 31 & 129 & 13 & $127^x + 32^y=(13z)^2$\\ [1ex] 
 \hline
\end{tabular}
 \caption{Some of the unsolvable cases of the equation $M_p^{x}  + (M_q+1)^y=(lz)^2$}
    \label{tab:my_label}
\end{table}

Now for possible extensions, the reader may try to solve the following Diophantine equations: \\
(i) $M_p^{x}+ (M_q + k)^y = z^2$, where $ k \geq 1$, and $M_p$ and $M_q$ are Mersenne primes. \\
(ii)  $M_p^{x}+ (M_q + 1)^y = z^n$, where $ n \geq 1$, and $M_p$ and $M_q$ are Mersenne primes and \\
(iii) $M_p^{x}+ (M_q + k)^y = z^n$, where $ k, n \geq 1$, and $M_p$ and $M_q$ are Mersenne primes.

\section{Acknowledgement}
 The author is grateful to Dr. Richa Sharma for her invaluable suggestions in preparing this article.

\end{document}